\documentclass[pdflatex,sn-nature]{sn-jnl}

\usepackage{graphicx}%
\usepackage{multirow}%
\usepackage{amsmath,amssymb,amsfonts}%
\usepackage{amsthm}%
\usepackage{mathrsfs}%
\usepackage{rotating}%
\usepackage[title]{appendix}%
\usepackage{xcolor}%
\usepackage{textcomp}%
\usepackage{etoolbox}%
\usepackage{manyfoot}%
\usepackage{booktabs}%
\usepackage{algorithm}%
\usepackage{algorithmicx}%
\usepackage{algpseudocode}%
\usepackage{listings}%
\usepackage{url}%
\usepackage{float}

\usepackage{geometry}%
\usepackage{hypcap}%
\usepackage{hyperref}%
\usepackage{breakurl}%
\usepackage{natbib}%
\usepackage{wrapfig}%

\usepackage{kotex}%
\usepackage{lineno}%
\theoremstyle{thmstyleone}%

\theoremstyle{thmstyletwo}%

\theoremstyle{thmstylethree}%

\DeclareMathOperator*{\argmin}{arg\,min}
\newcommand{\norm}[1]{\left\lVert#1\right\rVert}

\begin{document}

\title[Article Title]{
Cost-Effective Strategies for Infectious Diseases: A Multi-Objective Framework with an Interactive Dashboard
}

\author[1]{\fnm{Jongmin} \sur{Lee}}\email{ljm1729@konkuk.ac.kr}

\author[2]{\fnm{Renier} \sur{Mendoza}}\email{rmendoza@math.upd.edu.ph}

\author[2]{\fnm{Victoria May P.} \sur{Mendoza}}\email{vmpaguio@math.upd.edu.ph}

\author*[1]{\fnm{Eunok} \sur{Jung}}\email{junge@konkuk.ac.kr}

\affil[1]{\orgdiv{Department of Mathematics}, \orgname{Konkuk University}, \orgaddress{ \city{Seoul}, \postcode{05029}, \country{Korea}}}

\affil[2]{\orgdiv{Institute of Mathematics}, \orgname{University of the Philippines Diliman}, \orgaddress{\city{Quezon City}, \postcode{1101}, \country{Philippines}}}

\abstract{
During an infectious disease outbreak, policymakers face the critical challenge of balancing the healthcare burden with the socioeconomic consequences of social distancing. To facilitate this complex decision-making process, we developed an integrated framework that combines multi-objective optimization, economic evaluation, and an interactive dashboard. This dashboard allows stakeholders to input dynamic cost parameters and immediately derive cost-optimal intervention strategies. We validated this framework using data from the early COVID-19 outbreak in South Korea. Our results demonstrate that cost-optimal solutions remained remarkably consistent across a wide range of costs per infection (from 4,410 USD to 361,000 USD), exhibiting similar transmission reduction patterns. This suggests that the cost-optimal policy is relatively insensitive to variations in the estimated cost per infection, providing a robust basis for intervention planning. By synthesizing rapid optimization with rigorous economic evaluation, our framework offers a scalable tool for timely, evidence-based decision-making during the urgent situation of a pandemic.}




\maketitle

\section{Introduction}\label{sec1}

In recent decades, increasing air travel and the growing concentration of populations in urban areas have accelerated the global spread of infectious diseases, such as SARS, H1N1 pandemic influenza, MERS, and COVID-19 \cite{Bogoch2018,Baker2022}. Although traditional measures such as isolation, quarantine, and community containment were successfully implemented during the 2003 SARS outbreak \cite{WilderFreedman2020}, most countries were unable to suppress the spread of COVID-19 using these standard measures until the implementation of the vaccine \cite{LessonFromCOVID}. This discrepancy is largely attributed to the higher intrinsic transmissibility of COVID-19 and the unprecedented scale of global human mobility, which facilitated rapid cross-border seeding before governments recognized the spread of the disease. As an example, from the first confirmed case in the world, most countries had COVID-19 cases within 90 days, and unfortunately, there was no vaccine for the disease. Therefore, in the world, 1\% of the population was infected, and two million people died within one year. Consequently, there is an urgent need to plan an optimal non-pharmaceutical interventions (NPIs) strategy considering the country's situation \cite{Behavioural_environmental_social_systems_interventions}.

However, identifying an optimal strategy is challenging because multiple strategies may satisfy Pareto optimality, a state in which no further improvement can be made in one objective (e.g., reducing cumulative infections) without worsening another competing objective (e.g., minimizing intervention costs) \cite{TradeOff,HealthEconomy,EconoEvalDifficulty1,EconoEvalDifficulty2,EconoEvalDifficulty3}. In this multi-objective framework, a strategy is considered Pareto optimal if there is no other feasible strategy that performs better in at least one objective while remaining at least as good in all others. This approach allows us to identify the Pareto front, which represents the most efficient set of trade-offs from which policymakers can select based on societal priorities. Among the cost-optimal candidates, economic evaluation finds a cost-optimal solution. This strategy can be summarized by two sequential processes: (1) Finding the candidates that simultaneously minimize the number of infections and the economic impact of NPIs, and (2) determining a cost-optimal solution based on decision-makers' opinions. 

Existing research on infectious disease dynamics has extensively utilized compartmental, stochastic, agent-based, and network models to characterize transmission and evaluate intervention strategies \cite{compartmentModel1, compartmentModel2, stochasticModel1, stochasticModel2, agentModel1, agentModel2, networkModel1, networkModel2}. While traditional optimal control theory offers a robust framework for identifying intervention strategies, it often requires pre-defined objective weights that limit flexibility \cite{lenhart2007optimal,CBA2,optimizeIDM,OCEpi1,OCEpi2,OCEpi3}. To address this, multi-objective optimization has been adopted to derive Pareto-optimal strategies without the need for preset parameters or repetitive simulations \cite{MOOEpi1, MOOEpi2,MOOEpi3,MOOEpi4,MOOEpi5}. However, these theoretical results must be analyzed through economic evaluations to find a cost-optimal solution among the Pareto front \cite{CBA1, CBA3, CBA4, CBA5}. Despite these advancements, a critical gap remains: few studies simultaneously integrate multi-objective optimization with rigorous economic evaluation, and no study translates these complex outputs into accessible, real-time tools for policymakers \cite{MOODashboard1,MOODashboard2}. Our work addresses this issue by synthesizing modeling, optimization, and cost-benefit analysis into a single workflow, delivered via an intuitive dashboard designed for rapid, actionable decision-making during an epidemic.

This study introduces a research framework for addressing general infectious diseases. Subsequently, we present an application to the early phase of COVID-19 in Korea, which encompasses mathematical modeling, parameter estimation, multi-objective optimization, cost-benefit analysis, and the development of an interactive dashboard. This study makes three key contributions. First, we propose a novel methodological approach that applies multi-objective optimization to address the inherent trade-offs in policy interventions, a strategy that can be generalized to other infectious diseases beyond COVID-19. Second, we reveal that only a few cost-optimal patterns of transmission reduction exist in non-pharmaceutical intervention (NPI) strategies, depending on infection cost. Third, the interactive dashboard serves as an intuitive decision-support tool that provides optimal trajectories of transmission reduction over time. Using the multi-objective optimization framework, users can evaluate how varying NPI intensities influence both public health and economic outcomes. The remainder of the manuscript is organized into several sections. Section 2 outlines the methods of this study. Section 3 presents the results. Finally, Section 4 discusses the results, along with their implications and limitations, in comparison to existing literature. 

\section{Methods}\label{sec4}

\subsection{Workflow for this research}

\begin{figure}[H]
\centering
\includegraphics[width=\linewidth]{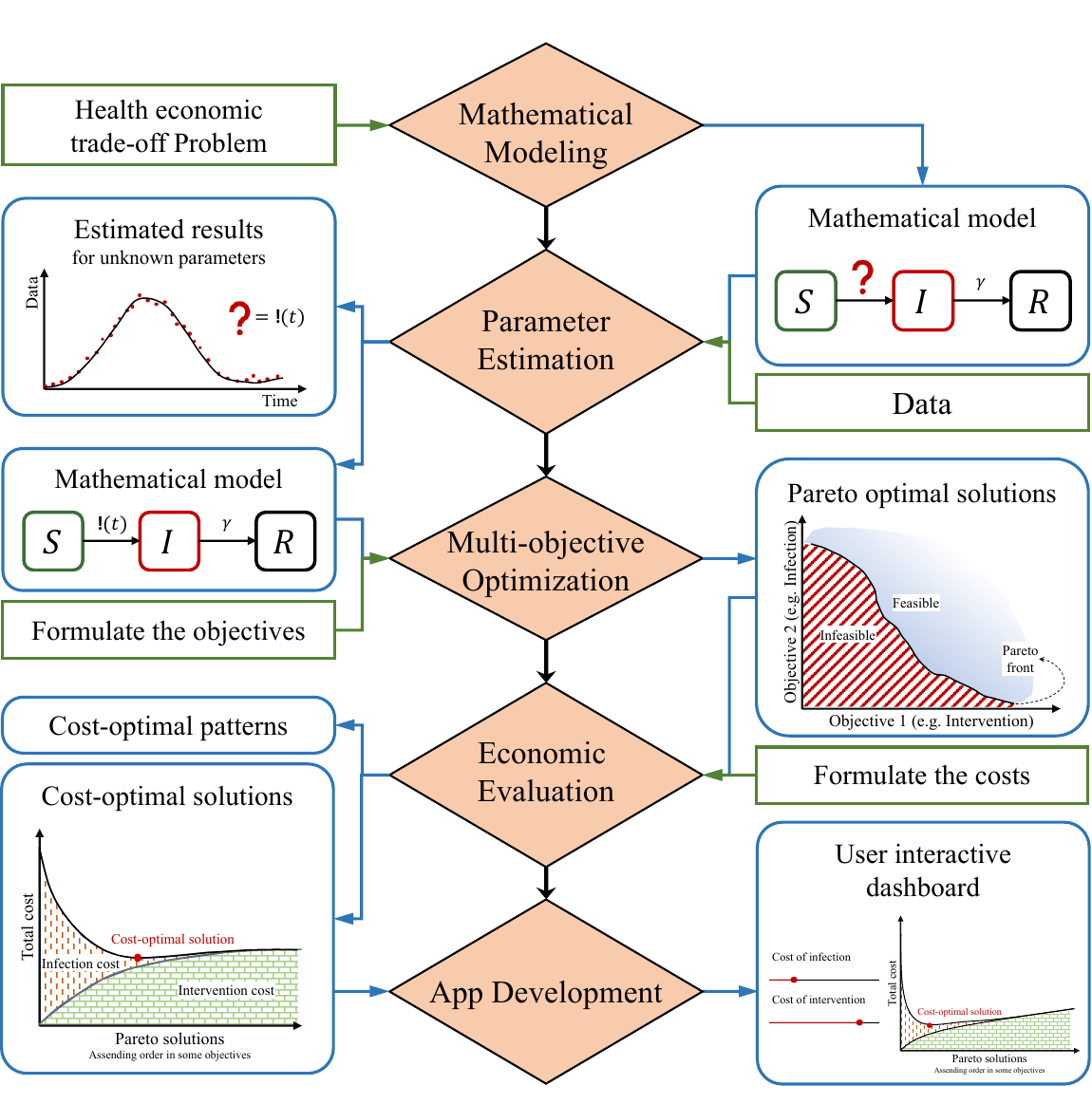} 
\caption{\textbf{Schematic overview of the integrated decision-support framework for pandemic management.} The framework consists of five sequential stages (orange diamonds) designed to suggest a tool for policy decision-making under a pandemic. 1) Mathematical Modeling for the target outbreak. 2) Parameter Estimation to estimate unknown parameters against observable data. 3) Multi-objective Optimization to derive the Pareto optimal solutions. 4) Economic evaluations to identify cost-minimizing intervention patterns among Pareto optimals. and 5) App Development resulting in a user-interactive dashboard for real-time strategy visualization under customized cost-related parameters.}\label{fig:01}
\end{figure}

To support decision-making in the early phases of infectious disease outbreaks, we developed a systematic research framework consisting of five sequential processes, as illustrated by the orange diamonds in Figure~\ref{fig:01}. In this study, we employed an augmented SIR-type compartmental model incorporating key epidemiological features, including the latent period, imported cases, isolation, and disease-induced mortality, specifically adapted to the early COVID-19 pandemic in South Korea. Using cumulative confirmed case data, we estimated key parameters, the average daily imported cases, the reduction in the infectious period due to NPIs, and the time-varying transmission reduction.

In the third step, multi-objective optimization is employed to identify the Pareto optimal solutions that simultaneously minimize the total number of infections and the average transmission reduction over time. To derive a cost-optimal solution from this Pareto set during a pandemic, a cost-benefit analysis is conducted using country-specific economic parameters. The framework culminates in the development of a decision-support tool (app development), resulting in a user-interactive dashboard. This dashboard enables end-users to adjust cost assumptions in real-time and visualize the corresponding cost-optimal strategies, providing a flexible and data-driven tool for policy decisions under various pandemic scenarios.

\subsection{Data}
The time-series data for confirmed cases and deaths in South Korea were obtained from Our World in Data \cite{owid-who-covid} and are available on our GitHub repository (\href{https://github.com/ljm1729-scholar/MOO_framework.git}{https://github.com/ljm1729-scholar/MOO\_framework}), along with the source code for our framework. The dataset covers the period starting from January 19, 2020, when the first case was officially reported.


\subsection{Mathematical modeling} \label{sec4_01} 

\begin{figure}[H]
\centering
\includegraphics[width=\textwidth]{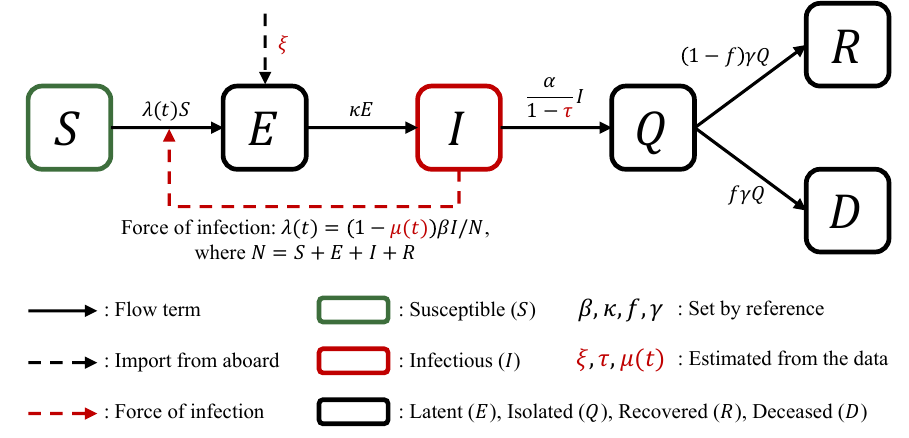}
\caption{\textbf{Mathematical model for infectious diseases.} 
The squares represent compartments of the mathematical model, and the black arrows represent flows between the compartments. The black dashed arrow represents external importation, which serves as a trigger for an epidemic. The red dashed arrow represents the force of infection, which drives the spread of the disease in a country. The parameters in black and red are the disease-related and estimated parameters, respectively.}\label{fig:02}
\end{figure}

Figure~\ref{fig:02} illustrates the susceptible ($S$), exposed ($E$), infectious ($I$), isolated ($Q$), recovered ($R$), and deceased ($D$) (SEIQRD) model used to investigate the early phases of COVID-19 in Korea. Note that $E$ is exposed to the disease and is infected; however, they are not infectious. $Q$ is isolated because they have been confirmed as infected cases by governmental or medical staff. Five disease-related parameters exist ($R_0$, $\kappa$, $\alpha$, $f$, and $\gamma$), where $R_0$ is the basic reproductive number of the disease, $1/\kappa$ is the average latent period, $1/\alpha$ is the average infectious period, $f$ is the fatality rate, and $1/\gamma$ is the average isolation period. Three policy-related parameters ($\mu(t)$, $\xi$, and $\tau$) exist. $\mu(t)$ represents the transmission reduction by NPIs, $\xi$ indicates the average number of unrecognized imported cases per day, and $\tau$ denotes the average reduction in the infectious period resulting from public health interventions. 

    The governing equations~(\ref{eq01}--\ref{eq06}) represent the changes in each compartment.
    \begin{align}
        \frac{dS}{dt} &=-\lambda (t) S, \label{eq01}\\
        \frac{dE}{dt} &= \lambda (t) S - \kappa E +\xi, \label{eq02} \\
        \frac{dI}{dt} &= \kappa E - \frac{\alpha}{1-\tau} I, \label{eq03} \\
        \frac{dQ}{dt} &= \frac{\alpha}{1-\tau} I - \gamma Q, \label{eq04} \\
        \frac{dR}{dt} &= (1-f) \gamma Q,\label{eq05} \\
        \frac{dD}{dt} &= f \gamma Q,\label{eq06}
    \end{align}
    where $\lambda(t)=(1-\mu(t))\mathcal{R}_0\frac{\alpha}{1-\tau} I/N$ is the force of infection, $\mathcal{R}_0$ is the basic reproduction number, $1/\alpha$ is the average infectious period, $\tau$ represents the average infectious period reduction that accounts for the shortened infectious period resulting from public health interventions, e.g. testing and contact tracing. $\mu(t)$ is the time-dependent transmission reduction resulting from policy interventions, and $N = S+E+I+R$. The parameter $\mu(t)$ consists of a set of values $\mu_i$, each representing the level of transmission reduction during a specific period. The index $i$ corresponds to the period starting on day $14 \times (i - 1)$ from the beginning of the simulation, with each period spanning two weeks.
    The details regarding $\mu(t)$ are provided in Appendix~\ref{asec1}. Notably, the infectious period does not begin with the onset of symptoms but rather with the start of infectiousness. For example, patients infected with COVID-19 can be infectious two days before developing symptoms \cite{beforeSymptom}. Figure~\ref{fig:02} presents a flowchart of the model.

\begin{table}[ht]
    \centering
    \caption{Parameter table for the SEIQRD mathematical model. The symbols, definitions, and values are presented in the table with corresponding references. $\mu(t)$, $\xi$, and $\tau$ are estimated from the cumulative confirmed data.}
    \begin{tabular}{clll}
        
        \toprule
        \textbf{Symbol} & \textbf{Definition} & \textbf{Value} & \textbf{References} \\
        \midrule 
        $R_0$  & Basic reproductive number  & 2.87  & \cite{R0COVID} \\
        1/$\kappa$ & Average latent period & 4 (day) & \cite{latentCOVID} \\
        1/$\alpha$ & Average infectious period  & 10 (day)  & \cite{infectiousCOVID} \\
        1/$\gamma$ & Average isolation period & 14 (day) & \cite{OptimalQuarantine1,OptimalQuarantine2} \\
        $f$ & Case fatality ratio & 0.0173 & \cite{owid-who-covid} \\
        $\mu(t)$ & Transmission reduction by NPIs over time & [0$\sim$0.95]*  & estimated \\
        $\xi$ & Average unrecognized imported cases per day & 0.2780* (person/day) & estimated \\
        $\tau$ & Average infectious period reduction by NPIs& 0.6218* & estimated \\
        \bottomrule
    \end{tabular}
    *: estimated parameter
    \label{tab:01}
\end{table}
    
    We utilized two global optimizers, the improved multi-operator differential evolution (IMODE) and Markov chain Monte Carlo (MCMC), to estimate the parameters \cite{IMODE, MCMC_package, MCMC}. IMODE ranked first in the congress on evolutionary computation (CEC) 2020 competition on bounded single-objective optimization algorithms. It combines the advantages of global and local search strategies, focusing on exploration at the beginning and exploitation at the end of the optimization process. IMODE excels at determining the optimal solutions for given function evaluations within a specified domain, without a preset initial point. We employed IMODE to minimize the difference between the model simulation results and the cumulative confirmed case data ($C_{data}(t)$). 
    \begin{equation}
    \min_{\theta}\sum_{t_0}^{t_{end}}\norm{\int_{t_0}^{t}\frac{\alpha}{1-\tau}I(s;\theta)ds - C_{data}(t)}_2 dt \label{eq07}
    \end{equation}
    where $\theta\in{\{\xi,\tau,\mu(t)\}}$ are policy-related parameters. The $\frac{\alpha}{1-\tau}I(t)$ represent the daily confirmed cases, which is the flow from $I$ to $Q$. The $\int_{t_0}^{t}\frac{\alpha}{1-\tau}I(s)ds$ represents the cumulative confirmed cases.

\subsection{Multi-objective optimization}\label{sec4_02} 
    We assume that interventions affecting the infection period and the number of imported cases remain fixed during an epidemic, as these are largely determined by a country’s healthcare infrastructure. In multi-objective optimization, we focus exclusively on changes in the transmission reduction parameter $\mu(t)$, which is influenced by government policies such as mask-wearing, social distancing, and gathering restrictions. The goal of optimization is to simultaneously minimize infections and the impact of transmission-related interventions by adjusting $\mu(t)$. Here, we minimize the number of infections instead of the number of confirmed cases used in the parameter estimation process because of the time delay between infection and confirmation. The expression $\int_{t_0}^{t_f}\lambda(t)S(t)dt$ reflects the total number of infections over the simulation period; the two objectives reflect intervention stringency and epidemic size. However, directly quantifying the relationship between $\mu(t)$ and the intervention cost is difficult. To address this, we assume that the intervention cost is proportional to its effectiveness and the stringency of the policy. Similarly, the infection cost is proportional to the total number of infections. Thus, the optimization aims to minimize both cost functions simultaneously. 
    $f_1(\mu(t))$ and $f_2(\mu(t))$ are represented by 
    \begin{equation}
        \argmin_{\mu(t)}
        \begin{bmatrix}
            f_1\left(\mu\left(t\right)\right) \\
            f_2(\mu(t)) 
        \end{bmatrix}
        =\argmin_{\mu(t)}
        \begin{bmatrix}
            \frac{\int_{t_0}^{t_f}\mu(t)dt}{t_f -t_0} \\
            \int_{t_0}^{t_f}\mu(t)\mathcal R_0\alpha I(t)\frac{S(t)}{N(t)}dt 
        \end{bmatrix}. \label{eq08}
    \end{equation}
    Note that $f_1$ and $f_2$ are proportional to the monetary cost but are not exact representations. The multi-objective optimization problem simultaneously minimizes $f_1$ and $f_2$, subject to the governing equations given by equations \eqref{eq01}--\eqref{eq06}.
    
    Multi-objective optimization identifies solutions near the Pareto curve; that is, a set of Pareto solutions in the objective plane composed of the co-domains of $f_1$ and $f_2$. A solution is considered Pareto optimal if it is not dominated by any other solution; that is, no other solution performs better in at least one objective without performing worse in another. We obtained the Pareto curve using the built-in function \textit{multiobjga}, which employs the non-dominated sorting genetic algorithm II (NSGA-II) in MATLAB \cite{multiobjga}. To ensure accuracy, we independently ran \textit{multiobjga} one thousand times and assembled the Pareto solution by locating the Pareto front. 

\subsection{Cost-benefit analysis}\label{sec4_05} 
    To determine a cost-optimal solution on the Pareto curve, a cost-benefit analysis is conducted by converting epidemiological outcomes into monetary terms. The total societal cost, $C_{\textrm{tot}}$, is defined as the sum of the intervention cost ($C_{\textrm{int}}$) and the infection cost ($C_{\textrm{inf}}$).
\subsubsection{Intervention cost $C_{\textrm{int}}$ formulation}
    The intervention cost is assumed to be proportional to the effectiveness and stringency of the transmission reduction strategy ($\mu(t)$). It is calculated by multiplying the maximum daily intervention cost, $C_1$, by the average reduction in relative intervention cost, $f_1 (\mu(t)$).
\begin{equation}
    C_{\textrm{int}} (\mu(t))  = C_1 f_1 (\mu(t))  =C_1  \left(\int_{t_0}^{f_f}\mu(t)dt/(t_f-t_0 )\right) \label{eq09}
\end{equation}
The parameter $C_1$ is derived from the country's Gross Domestic Product (GDP) and the maximum projected GDP reduction  ($\textrm{GDP}_{\textrm{MaxRed}}$) \cite{GDP, GDPforecast} due to non-pharmaceutical interventions.
\begin{equation}
    C_1=\textrm{GDP}\times\textrm{GDP}_{\textrm{MaxRed}}.\label{eq10}
\end{equation}
This $\textrm{GDP}_{\textrm{MaxRed}}$ is derived from the difference between the actual GDP during the COVID-19 pandemic and the GDP of South Korea in 2020, projected by the OECD in 2019. The calculation is proposed in the Supplementary.

\subsubsection{Infection cost $C_{\textrm{inf}}$ formulation}
The infection cost is proportional to the total number of infections and is formulated as follows:
\begin{equation}
    C_{\textrm{inf}} (\mu(t))  = C_2 f_2 (\mu(t))= C_2 \int_{t_0}^{t_f}\lambda(t)S(t)dt.\label{eq11}
\end{equation}
Here, $C_2$ represents the cost per infection, which incorporates both the average hospitalization cost ($C_{\textrm{H}}$) and the economic cost of mortality, calculated using the case fatality rate ($f$) and the value of a statistical life (VSL) \cite{CBA3,CBA4,CBA5,CBA_HospitalizationCost,CBA_ValueOfStatisticalLife}:
\begin{equation}
    C_2  = C_\textrm{H}+f \times \textrm{VSL}.\label{eq12}
\end{equation}

\subsection{App development}\label{sec4_06}
    Assigning fixed values to the parameters in equations \eqref{eq10} and \eqref{eq12} is challenging. For example, $\textrm{VSL}$ can change according to average age, wage, income elasticities, and ethical considerations \cite{VSL_hard1, VSL_hard2, VSL_lowmiddle, VSL_hard3, VSL_world1}. To accommodate this variability, we developed a dashboard that enables users to select their own cost-related parameters and obtain cost-optimal results within seconds. The web dashboard, built using the Shiny package in Python, includes both a mathematical model simulator and \href{https://jongmin-lee.shinyapps.io/demomookorea/}{cost-optimal intervention policy simulator}. In the mathematical model simulator, users can alter various parameters of the SEIQRD model and instantly view the simulation results. In the policy simulator, users can adjust five economy-related parameters that serve as inputs to $C_1$ and $C_2$. This real-time optimization simulator is feasible because multi-objective optimization does not require a weight parameter for each objective.

\section{Results}\label{sec2}

\subsection{Parameter estimation} \label{sec2_02} 
Table~\ref{tab:01} lists the parameters of the mathematical model. Parameters $R_0$, $1/\kappa$, $1/\alpha$, and $1/\gamma$ are epidemiological quantities that characterize the disease and whose values can be obtained from relevant references. However, $\mu(t)$, $\xi$, and $\tau$ may vary across countries or change over time, depending on the policies or interventions implemented during a given period. Therefore, these unknown parameters must be estimated. We utilized a hybrid parameter estimation scheme using two global optimizers to obtain a posterior distribution: the improved multi-operator differential evolution (IMODE) and the Metropolis-Hastings (MH) algorithm \cite{IMODE}, which is a Markov chain Monte Carlo (MCMC) method \cite{MCMC}. Table~\ref{tab:01} lists the estimated values of the policy-related parameters obtained by fitting the model to the cumulative confirmed case data. The average number of imported cases per day ($\xi$) was 0.2780, which is approximately one person every four days. The infectious period reduction ($\tau$) was 62.18\%, indicating that the contact tracing or testing policy reduced the infectious period by this value. The transmission reduction parameter, which has a value between 0 and 0.95, is estimated every two weeks as government transmission reduction policies changed frequently. Note that $\mu$ is set to zero during the first two weeks because confirmed cases had not yet been detected. Appendix~\ref{asec2} presents the correlations between the estimated parameters derived from the MCMC chain. 

\begin{figure}[H]
    \centering
    \includegraphics[width=\textwidth]{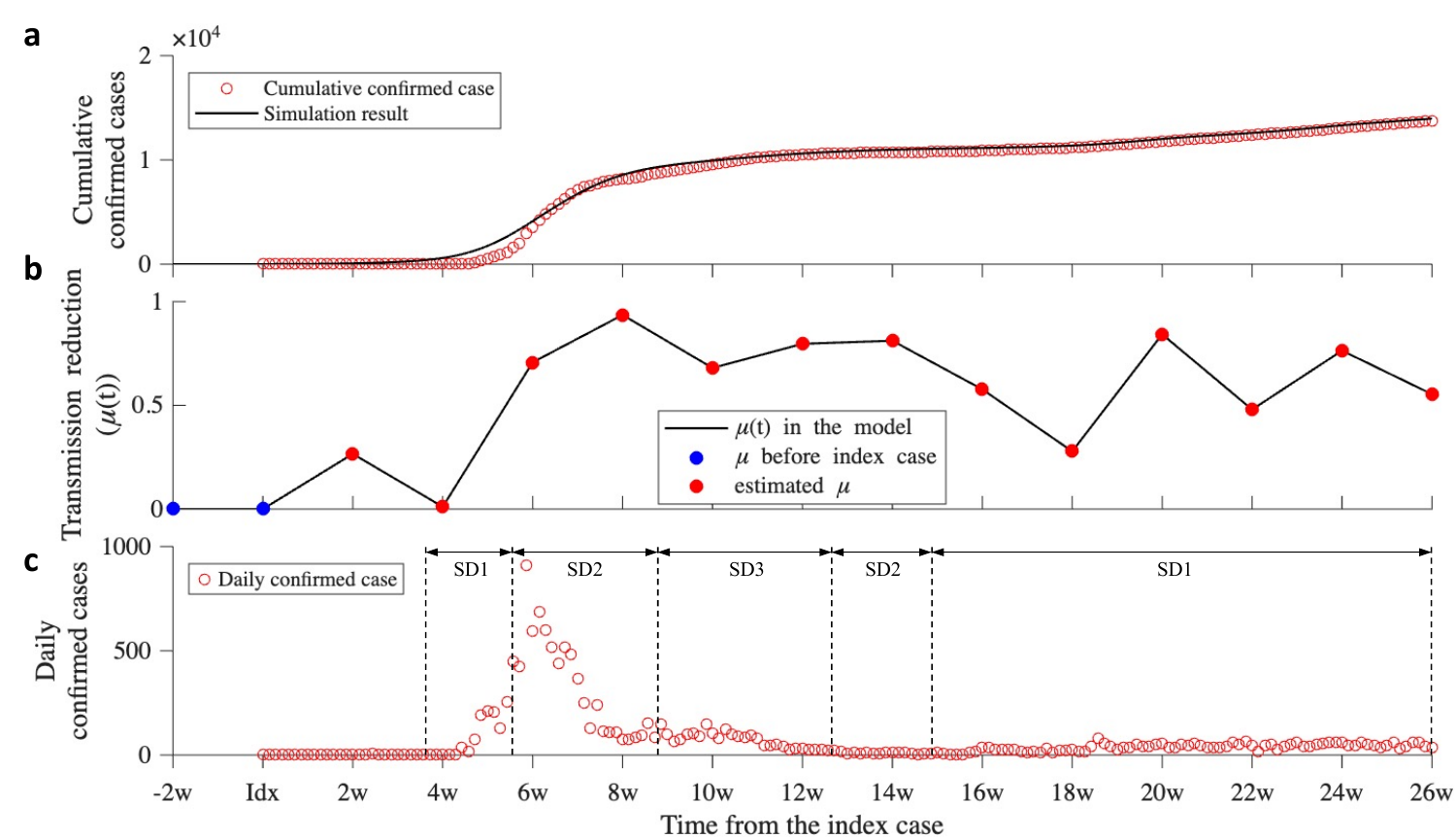}
    \caption{\textbf{Data-fitting results.} 
    (a) Cumulative confirmed cases and simulation results. (b) Estimated values of $\mu(t)$ over time. (c) Daily confirmed cases and corresponding implemented social distancing strategy.
    }\label{fig:03}
\end{figure}

The visual correspondence observed in our data-fitting result confirms that fluctuations in transmission rates, $\mu(t)$, reflect dynamic changes in population behavior driven by social distancing (SD) strategies, thereby providing a data-driven relationship between transmission reduction and SD levels. This relationship is based on the model’s ability to capture the epidemiological trends in cumulative COVID-19 cases in South Korea over a 26-week period (Figure~\ref{fig:03}a). Central to this calibration is the estimated time-varying transmission coefficient, $\mu(t)$, which reveals a sharp, step-like increase from a near-zero base to its peak during the 5th and 6th weeks, directly correlating with the observed exponential surge in cases (Figure~\ref{fig:03}b). A temporal analysis of daily confirmed cases alongside government-implemented social distancing (SD) levels further elucidates these patterns; the primary epidemic wave peaked during the less restrictive SD Level1 (SD1) period, while the subsequent deceleration between weeks 7 and 14 strongly coincides with the transition to more stringent SD Level2 (SD2) and SD Level3 (SD3)  measures (Figure~\ref{fig:03}c), validating the model's alignment with structural policy shifts.

\subsection{Characteristics of optimization: Pareto curve}\label{sec2_03} 
\begin{figure}[H]
\centering
\includegraphics[width=.98\textwidth]{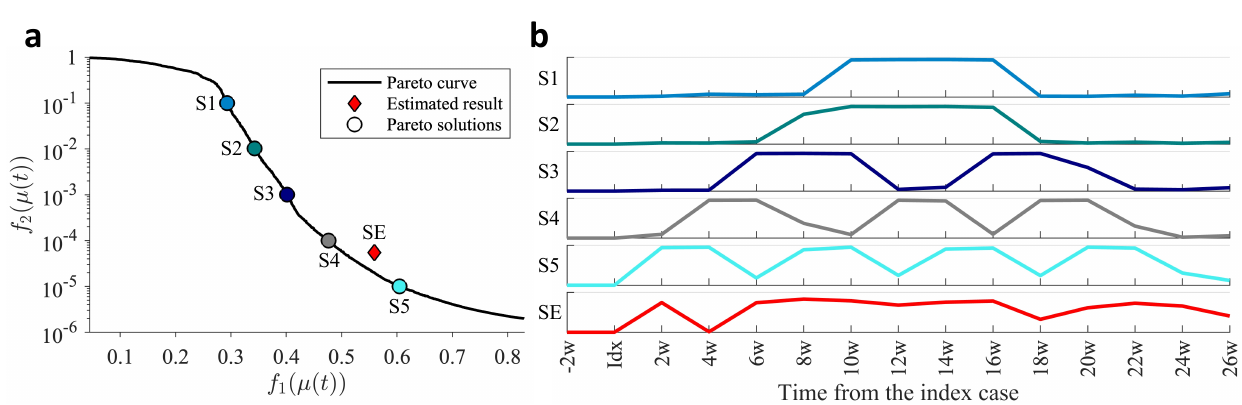}
\caption{\textbf{Pareto-optimal intervention strategies.} (a) Objective space and Pareto frontier: The black solid line represents the Pareto-optimal frontier, illustrating the trade-off between the average effectiveness of transmission reduction ($f_1\left( \mu \left(t\right)\right)$) and the proportion of the population infected. The colored circles (S1–S5) denote five representative Pareto-optimal solutions selected based on specific cumulative infection rate targets: 10\% (S1), 1\% (S2), 0.1\% (S3), 0.01\% (S4), and 0.001\% (S5). The red diamond (SE) indicates the estimated real-world strategy from South Korean data, which sits above the Pareto curve, reflecting the inherent discrepancy between empirical implementation and theoretical optimality. (b) Temporal transmission-reduction profiles: Time-series representation of the specific intervention patterns corresponding to the points in panel (a). The profiles show that achieving more stringent public health targets (moving from S1 to S5) requires initiating high-intensity interventions earlier in the outbreak. The red curve (SE) displays the historical transmission reduction estimated during the 26-week study period.}\label{fig:04}
\end{figure}

Figure~\ref{fig:04}a reveals a clear trade-off between public health preservation and intervention burden, demonstrating that incremental improvements in infection control require progressively higher levels of transmission reduction. The required average transmission reduction increases from 32.13\% (S1) to 58.32\% (S5) as the target infection threshold is tightened from 10\% to 0.001\% of the population, respectively. While the estimated strategy from empirical data (SE) yielded a 0.0279\% infection rate with a 53.05\% average transmission reduction, it is positioned within the feasible region but off the Pareto curve. This gap indicates that the historical response, while effective, did not reach the theoretical frontier of optimality, suggesting that more efficient trade-offs could have been achieved through optimized resource allocation and timing.

Temporal analysis of the optimal strategies in Figure~\ref{fig:04}b indicates that the timing of intervention onset is the critical factor in achieving stringent public health targets. Specifically, while strategies with higher infection tolerances (S1) allow for delayed responses, the most stringent strategy (S5) requires the initiation of high-intensity transmission reduction as early as the 2nd week after the index case is detected. The estimated real-world strategy (SE) have same delay with S3 strategy although the closest point in Figure~\ref{fig:04}a is S4. The model suggests that achieving these targets does not necessitate an unbroken, maximum-intensity intervention; rather, it favors adaptive fluctuations in policy stringency throughout the target period to minimize the total intervention burden while keeping infections within the prescribed limits.

\subsection{Cost-optimal solutions} \label{sec2_04}
\begin{figure}[H]
\centering
\includegraphics[width=.98\textwidth]{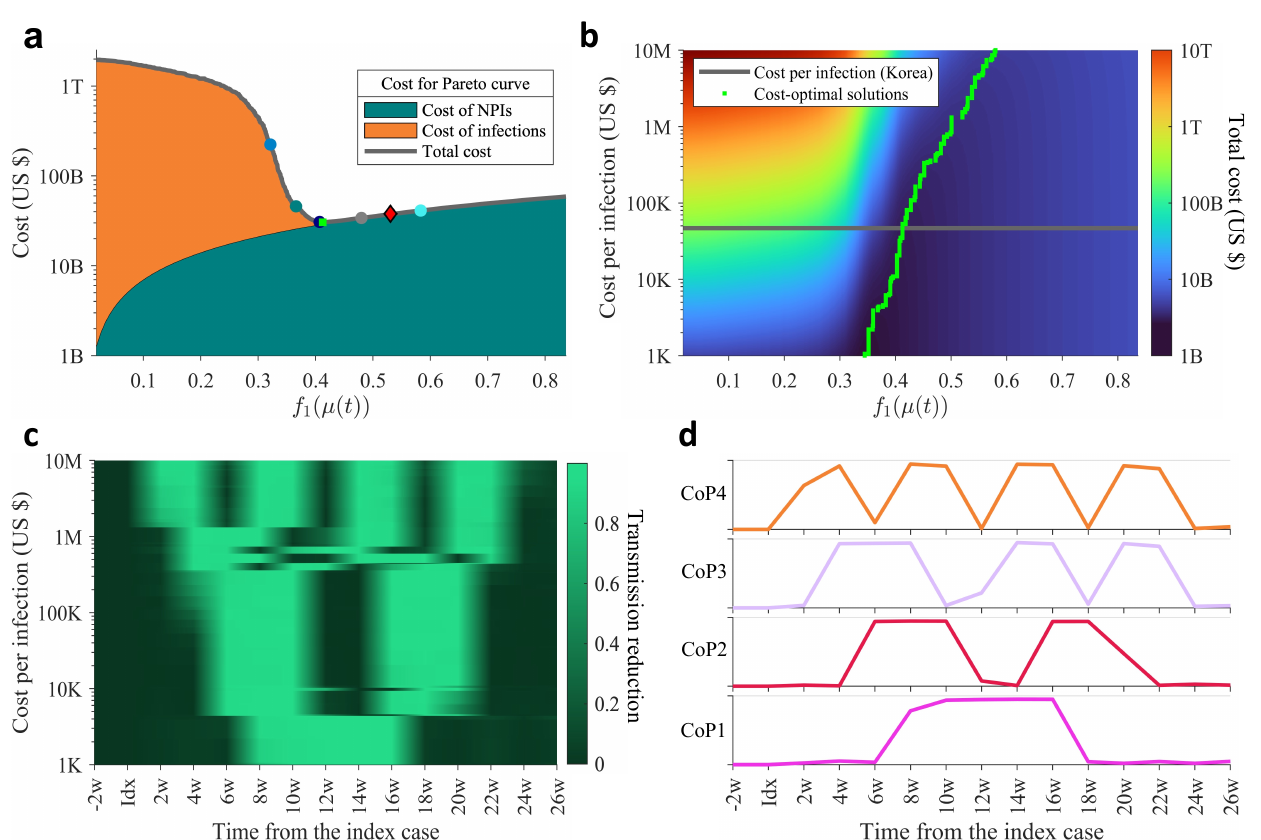}
\caption{\textbf{Cost-benefit analysis results with Pareto solution.} (a) Total costs along the Pareto curve under the cost-benefit analysis. The orange area, green area, and gray line represent the infection cost, transmission reduction-related intervention cost, and total cost, respectively. (b) Total costs along the Pareto curve for different values of cost per infection. The green points represent the cost-optimal solutions for varying infection costs. The gray line corresponds to panel (a). (c) Cost-optimal solutions for the cost per infection. The line parallel to the x-axis represents the cost-optimal solution for each infection cost and corresponds to the green points in the panel (b). (d) Cost-optimal policies for each cost per infection range.}\label{fig:05}
\end{figure}

The cost-benefit analysis is conducted to identify a cost-optimal solution among the Pareto-optimal frontier, which minimizes the total cost. Figure~\ref{fig:05}a demonstrates this relationship by decomposing the total societal cost ($C_{\textrm{tot}}$) into its two constituent parts: the cost of infections ($C_{inf}$) and the cost of non-pharmaceutical interventions ($C_{int}$). For the specific economic context of South Korea with a cost per infection of 39,213 USD, the cost-optimal solution (represented by a green square) corresponds to an average transmission reduction of 41.25\%, resulting in a total cost of 30.6 billion USD. Comparatively, the estimated historical strategy (SE) resulted in a higher total cost of 37.9 billion USD, representing a 7.3 billion USD gap from the theoretical optimum.

A sensitivity analysis across a wide economic range reveals that the cost-optimal transmission reduction level is highly dependent on the assigned cost per infection. As shown in the heatmap in Figure~\ref{fig:05}b, the optimal $f_1(\mu(t))$ shifts rightward as the cost per infection increases from 1,000 USD to 10 million USD, indicating that higher valuations of health outcomes necessitate more stringent intervention levels to minimize total societal costs.

Despite the existence of over a thousand Pareto-optimal solutions, fewer than 100 are identified as cost-optimal across the range of infection costs. These cost-optimal exhibit a remarkably high degree of structural similarity within specific cost intervals. For example, when the cost per infection ranges from 4.41K USD to 361K USD, the solutions exhibit a remarkably similar strategic structure. This high similarity in a range of cost per infection suggests that a robust policy remains optimal regardless of minor economic fluctuations, leading us to define the Cost-Optimal Pattern (CoP). These CoPs characterize stable intervention behaviors, such as the optimal timing of initiation and duration of stringency, providing a practical reference for decision-makers.

As illustrated in Figure~\ref{fig:05}c and \ref{fig:05}d, as the cost per infection rises, the optimal strategy transitions from a delayed, single-peak intervention (CoP1) to a pattern of multiple, early-onset high-intensity periods (CoP4). Specifically, for the lower cost range (CoP1), the optimal strategy begins with strong interventions at the 6th week, whereas for the highest cost range (CoP4), the strategy initiates strong measures immediately at week 0 and maintains a periodic, high-stringency profile to prevent any significant surge in cases.

\begin{table}[ht]
  \centering
  \caption{\textbf{Cost-optimal pattern (CoP) from the index case for cost per infection.} }
  \label{tab:02}
      \begin{tabular}{@{} l l c c c c l l @{} }
        \toprule
        Cost per  &
        \multirow{2}{*}{Scenario} &
        \multicolumn{4}{c}{Strategy pattern (week)} &
        Total cost &
        Total \\
        \cmidrule(lr){3-6}
        infection (USD)& & Begin & Increase & Strong & Decrease & (USD) & infection \\
        \midrule
        \multirow{2}{*}{1.33M --}    & \multirow{2}{*}{CoP4} & \multirow{2}{*}{0} & 0--1, 6--7,  & 2--3,8--9, & 4--5,10--11, & 38.96B --      & 0.00103\% --  \\
         & & & 12--13,18--19 & 14--15,20--21 & 16--17,22--23 & 46.12& 0.00338\%\\
         
        \multirow{2}{*}{361K -- 1.33M}   & \multirow{2}{*}{CoP3} & \multirow{2}{*}{2} & 2--3, 10--13,          & 4--7, 14--15,    & 8--9, 16--17    & 38.96B --      & 0.00338\% --   \\
         & & & 18--19 & 20--21& 22--23 & 37.29B & 0.00678\% \\
         
        \multirow{2}{*}{4.41K -- 361K}       & \multirow{2}{*}{CoP2} & \multirow{2}{*}{4} & \multirow{2}{*}{4--5, 14--15}          & \multirow{2}{*}{6--9, 16--17}    & \multirow{2}{*}{10--11, 18--21}          & 28.04B --     & 0.00678\% --      \\
         & & & & & & 37.29B  &0.467\%\\
         
         \multirow{2}{*}{-- 4.41K}       & \multirow{2}{*}{CoP1} & \multirow{2}{*}{6} & \multirow{2}{*}{6-7} & \multirow{2}{*}{8--15}    & \multirow{2}{*}{16--17}          & 25.62B --     &0.467\% --     \\
         & & & & & & 28.04B  &2.62\%\\
        \bottomrule
      \end{tabular}
    \end{table}

Table~\ref{tab:02} presents the CoPs corresponding to the different ranges of costs per infection. Solutions within each CoP exhibit similar characteristics in terms of transmission reduction strategies, specifically, the number of intervention periods, duration of strong interventions, and timing of their initiation and termination. This table may serve as a practical reference for decision-makers in selecting a cost-optimal strategy based on their assessments of infection and intervention costs.

\subsection{User-interactive dashboard} \label{sec2_05}


As intervention and infection costs are not only difficult to estimate, but also vary widely across settings and individuals, we developed an \href{https://jongmin-lee.shinyapps.io/demomookorea/}{interactive dashboard} using Python Shiny, which allows users to adjust cost-related parameters, including GDP, GDP reduction, the value of a statistical life (VSL), and the fatality rate. The dashboard is published on two different servers: one is \href{https://jongmin-lee.shinyapps.io/demomookorea/}{Posit server, which is provided in the Shiny package by default} and the other is \href{http://www.infdissim.com}{the server provided by the authors of this paper}. 


\section{Discussion}\label{sec3}

    We formulated a mathematical model that incorporates the importation of infected individuals from abroad, transmission reduction, and infectious period reduction, which are key factors influenced by intervention policies during a pandemic. As the effects of interventions are initially unknown, these relevant parameters were estimated using a hybrid parameter estimation method that combines machine learning-based global optimization \cite{IMODE} and a statistical optimization method \cite{MCMC,MCMC_package}. By utilizing a mathematical model with estimated parameters, multi-objective optimization \cite{multiobjga} obtains Pareto solutions for transmission reduction policies that simultaneously minimize intervention and infection costs. Thereafter, a cost-benefit analysis determines the cost-optimal solution among the Pareto solutions based on the cost per infection; for example, a higher cost per infection necessitates a stronger cost-optimal intervention. To accommodate variability in cost assumptions, we developed a web-based dashboard that enables users to customize the cost, thereby obtaining the corresponding cost-optimal strategies. 

    The results obtained from the case study of South Korea can be interpreted through both mathematical and biological perspectives. Mathematically, the Pareto front represents the possible solutions that minimize the impact of infection and intervention at the same time (Figure~\ref{fig:04}a). Moreover, reducing the transmission from the cost-optimal solution makes more economic burden than high-intensity policies (Figure~\ref{fig:05}a). Biologically, the identified strategies (S1–S5) emphazise that early-phase control is governed by the need to keep $R_t$ below 1; once the epidemic surpasses a critical threshold, the cost of returning to a manageable state increases rapidly. This finding is consistent with previous studies that emphasize the importance of timing in NPI implementation \cite{NPIsEffect1,NPIsEffect2}. Specifically, our cost-optimal results for South Korea, which suggest a transition toward stronger interventions as the valuation of health outcomes (infection cost, $C_\textrm{inf}$) rises, align with the suppression versus mitigation trade-offs documented in early 2020 modeling efforts \cite{stochasticModel1,chin2021effect,sharma2021understanding}. While much of the existing literature provides a static, retrospective scoreboard, our framework advances these findings by offering a dynamic, proactive tool that bridges the gap between theoretical modeling and actionable policy.

    The elements of the framework can be modified according to different situations of the pandemic. First of all, the mathematical model should be adjusted if the transmission dynamics, fatality rates, pharmaceutical interventions, age structure, or other characteristics of a pandemic need to be considered. For example, when applying the framework to an emerging or unknown infectious disease, informed assumptions regarding the epidemiological characteristics of the disease and country-specific capabilities must be made. Based on these assumptions, the framework can simulate the introduction of a disease into a population via importation and identify Pareto-optimal and cost-optimal intervention strategies. To obtain results for different infectious diseases, users can modify the disease-related parameters using values from existing literature \cite{ReproductionNumberRID,IncubationPeriodRID,SerialIntervalRID} and then apply the framework accordingly. Similarly, to adapt the analysis to different countries, users can adjust the population- and policy-related parameters based on country-specific data. In this manner, cost-optimal solutions for other countries can be derived using their respective data.

    The other components of the framework can be adapted to align with specific research goals or data availability. For instance, the parameter estimation method \cite{DESS,RstanBayesian,ChowellBayesian} can be upgraded to more advanced Bayesian techniques or omitted entirely if all model parameters are pre-established. Regarding the optimization phase, the current objectives (cumulative infections and transmission reduction) can be substituted with alternative performance metrics, such as peak hospital demand or specific economic indicators, depending on the policymaker's priorities. This modularity ensures that the framework remains applicable across diverse public health scenarios and evolving policy requirements. The optimization method can also be modified to a state-of-the-art or appropriate algorithm for the specific problem \cite{RecentMOO1,RecentMOO2,RecentMOO3}. Cost-benefit analysis can be replaced with other economic evaluation methods by the objectives of researchers who want to use this framework. For example, framework users may modify the intervention costs, including contact tracing policies or pharmaceutical interventions, or use quality-adjusted life years or disability-adjusted life years to evaluate the impact on the economy.

    There are many formats available to provide the results, including Excel, raw code, installable program, a package, and a web dashboard. If your institute or country is familiar with a specific format, it is advisable to use that format for your own purposes. We publish a web-based dashboard using the Shiny package because it does not require the installation of any programs and is not restricted by the device. One challenge is the need for a reliable server to run our dashboard. A key feature of the published results is that users of the program can adjust the cost-related parameters and observe the corresponding results in real-time. Consequently, individuals who are not familiar with the details of the mathematical model or the optimization process can easily access and interpret the results. 

    While traditional optimal control theory provides a robust and efficient framework for identifying optimal intervention trajectories under a pre-defined single objective functional, it requires the fixed specification of weights between competing objectives prior to the optimization. In infectious disease modeling, where economic metrics such as the value of a statistical life (VSL) or the cost of treatment are often subject to high uncertainty and varying societal priorities, this requirement can be restrictive. To address this, we employ multi-objective optimization to identify the Pareto front, which maps the entire landscape of efficient trade-offs regardless of specific cost assumptions. By doing so, the cost-optimal solution for any given set of weights can be identified as a subset of the Pareto front. This enables a dynamic cost-benefit analysis across a broad range of economic scenarios (e.g., varying the cost per infection from 1K USD to 10M USD) without requiring repeated optimization runs for every parameter adjustment.

    In this study, we considered only transmission-related interventions to obtain Pareto-optimal solutions. However, optimal quarantine and testing strategies are equally important during pandemics \cite{OptimalQuarantine1,OptimalQuarantine2}. Although our proposed framework provides a method to analyze cost-free and cost-optimal intervention strategies, its practical implementation is limited. First, we adopted the established SEIQR model, which enables rapid experimentation and outcome generation. Although this framework emphasizes structural design, it remains flexible. Users can readily incorporate more complex models \cite{PracticalODE} without requiring algebraic derivations, owing to the implementation of the metaheuristic algorithms \cite{IMODE, multiobjga}. Second, different Pareto solutions may emerge if users want to use different mathematical models, objective functions, or certain parameters. However, once Pareto solutions are calculated, users can adjust only the cost-related parameters to obtain user-defined cost-optimal solutions. Third, the theoretical solution did not provide specific policies for achieving the suggested transmission reduction levels. Policymakers should refer to other studies exploring the relationship between transmission reductions and specific policies \cite{NPIsEffect1,NPIsEffect2,NPIsEffect3}. Finally, we considered only infection costs, excluding factors such as medical resources or potential overburdening, and detailed the costs of other interventions, such as limitations in gathering, quarantine, and testing. Nevertheless, users can modify the assumed costs of interventions or infections to incorporate these factors and obtain corresponding cost-optimal solutions.

    By operating complex epidemiological dynamics into a transparent, cost-effective optimization framework, our study provides decision-makers with an actionable control panel rather than a retrospective scoreboard. This integrative framework synthesizes epidemiological modeling, multi-objective optimization, and economic evaluation into a user-interactive tool that facilitates real-time exploration of what-if scenarios and cost assumptions. As a result, it transforms conceptual trade-offs into visible, evidence-based guidance that is interpretable by policymakers, public health practitioners, and clinical planners. Health authorities can identify intervention intensities that mitigate the disease overburden while preserving critical medical capacity; clinical stakeholders can be equipped with reasonable evidence for allocating resources—such as personal protective equipment, hospital beds, antibiotics, and vaccines—where they will yield the highest impact; and political decision-makers are empowered with transparent, data-driven justifications for public health measures that require sacrifices in civil liberties and economic considerations. Because all model assumptions and cost parameters are explicitly adjustable, stakeholders can test policies—before public trust evaporates, intensive care units are overwhelmed, or budget limits are reached—and converge on solutions that are politically viable based on theoretically optimal policy. In doing so, our framework advances pandemic governance from reactive crisis management to proactive, scenario-informed stewardship of public health, healthcare infrastructure, and societal resilience.

    \section*{Conclusion}

\backmatter

\bmhead{Supplementary information}

\section*{Declarations}
\subsection*{Funding}
This research was supported by the Bio\&Medical Technology Development Program of the National Research Foundation (NRF), funded by the Korean government (MSIT) (RS-2023-00227944). Additional support was provided by the Korea National Research Foundation (NRF) grant funded by the Korean government (MEST) (NRF-2021R1A2C100448711).

\subsection*{Data and material availability}
The confirmed cases and deaths obtained from Our World in Data \href{https://ourworldindata.org/coronavirus}{https://ourworldindata.org/coronavirus} \cite{owid-who-covid} have been uploaded to GitHub. The dashboard is available at \href{https://jongmin-lee.shinyapps.io/demomookorea/}{https://jongmin-lee.shinyapps.io/demomookorea/} or \href{http://www.infdissim.com}{http://www.infdissim.com}. 

\subsection*{Code availability}
The results can be generated using the code attached to GitHub 
\href{https://github.com/ljm1729-scholar/MOO_framework.git}{https://github.com/ljm1729-scholar/MOO\_framework}.

\subsection*{Author contribution}
\begin{itemize}
\item Conceptualization: Jongmin Lee, Renier Mendoza.
\item Methodology: Jongmin Lee, Renier Mendoza, Victoria May P. Mendoza.
\item Software and Validation: Jongmin Lee.
\item Formal Analysis and Investigation: Jongmin Lee, Renier Mendoza.
\item Data Curation: Jongmin Lee, Renier Mendoza, Victoria May P. Mendoza.
\item Visualization: Jongmin Lee.
\item Writing—Original Draft: All authors.
\item Supervision: Eunok Jung.
\item Project Administration: Eunok Jung.
\item Funding Acquisition: Eunok Jung.
\end{itemize}

\subsection*{Ethics declarations}
\subsubsection*{Competing interests}
The authors declare no competing interests.

\newpage
\begin{appendices}
\section{Transmission reduction $\mu(t)$}\label{asec1}
    The time-dependent transmission-reduction function $\mu(t)$ is defined using the linear interpolation of $\mu_i$, that is, the transmission reduction at time $t=i\times 14$: 
    \begin{equation}
        \mu(t)=\mu_{n+1} - \left(n-\frac{t}{14} \right) \left( \mu_{n+1}-\mu_{n}\right),
    \end{equation}
    where $n=\lceil\frac{t}{14}\rceil$ and $\mu_0=0$. We assume that transmission reduction changes every two weeks in accordance with government announcements regarding adjustments in intervention strategies. We further assume that the government implemented testing and control on imported cases. Parameter $\xi$ represents the mean number of imported cases per day. Imported cases are affected by screening measures at airports and borders. At the beginning of the simulation, we assumed no infected cases, as the outbreak was triggered by the imported cases. The parameter $\tau$ represents the infectious period reduction due to government testing and contact tracing efforts. We estimated these policy-related parameters using data on confirmed cases.
\section{Reproduction number}
The effective reproduction number, $R_t (t)$, is derived using the Next-Generation Matrix (NGM) method \cite{diekmann2010construction}. 
The appearance of new infections ($\mathcal{F}$) and the transfer of individuals out of these compartments ($\mathcal{V}$) are defined as follows:
\begin{equation*}
\mathcal{ F}=\begin{bmatrix}\lambda(t)S\\0\end{bmatrix},\mathcal{V}=\begin{bmatrix}\kappa E -\xi\\-\kappa E +\frac{\alpha}{1-\tau} I\end{bmatrix}
\end{equation*}
Substituting the force of infection $\lambda(t)=\left(1-\mu(t)\right) R_0  \frac{\alpha}{1-\tau}  \frac{I}{N}$:
\begin{equation*}
\mathcal{ F}=\begin{bmatrix}\left(1-\mu(t)\right) R_0  \frac{\alpha}{1-\tau}  \frac{I}{N}S\\0\end{bmatrix}
\end{equation*}
By computing the Jacobians $F$ and $V$ at the Disease-Free Equilibrium (DFE), where $S \approx N$, the Next-Generation Matrix $K = FV^{-1}$ is obtained as follows:
\begin{align*}
    F&=\frac{\partial \mathcal F}{\partial(E,I)}=\begin{bmatrix}0 & \left(1-\mu(t)\right) R_0  \frac{\alpha}{1-\tau}\\0 & 0\end{bmatrix} \\
V&=\frac{\partial \mathcal V}{\partial(E,I)}=\begin{bmatrix}\kappa & 0\\-\kappa & \frac{\alpha}{1-\tau}\end{bmatrix}
\end{align*}

The inverse of $V$ is:
\begin{align*}
    V^{-1}&=\begin{bmatrix}1/\kappa & 0\\1 & \frac{1-\tau}{\alpha}\end{bmatrix} \\
    K=FV^{-1}&=\begin{bmatrix}(1-\mu(t))R_0 & R_0\\0 & 0\end{bmatrix}.
\end{align*}

The effective reproduction number $R_t (t)$ is the spectral radius (the largest eigenvalue) of matrix $K$:
$R_t (t) =\rho(FV^{-1}) = (1-μ(t)) R_0$.
This result shows that the reproduction number is primarily driven by the basic reproductive number $R_0$ and scaled by the transmission reduction factor $(1-\mu(t))$. Note that the epidemic is suppressed ($R_t (t)<1$) when the intervention intensity $\mu(t)$ exceeds the threshold $1-1/R_0$ . For the case of COVID-19, the threshold value for $\mu(t)$ is estimated to be approximately 0.6516.

\section{Hybrid parameter estimation}\label{asec2}
As detailed in the Parameter Estimation section, IMODE is a heuristic global optimization algorithm designed to identify optimal solutions within a specified domain without a preset initial point. Because the IMODE algorithm effectively navigates the complex loss landscape to identify the most suitable parameter regions, its results provide a scientifically grounded basis for the subsequent statistical inference. IMODE is qualified by winning the 2020 Congress on Evolutionary Computation (CEC) competition, which requires solving 40 problems with a high-dimensional domain. Our optimization problem has 15 dimensions, which is fewer than the competition required. The optimized values from IMODE were utilized as the mean for the prior distributions in the MCMC process, ensuring that the Markov chains were initiated within high-probability regions of the posterior distribution. This hybrid approach mitigates the risk of the MCMC process becoming trapped in local modes and significantly enhances the efficiency and convergence of the parameter estimation.

\subsection{Improved multi-operator differential evolution}\label{asec2_1}
\begin{figure*}[h]
    \centering
    \includegraphics[width=\textwidth]{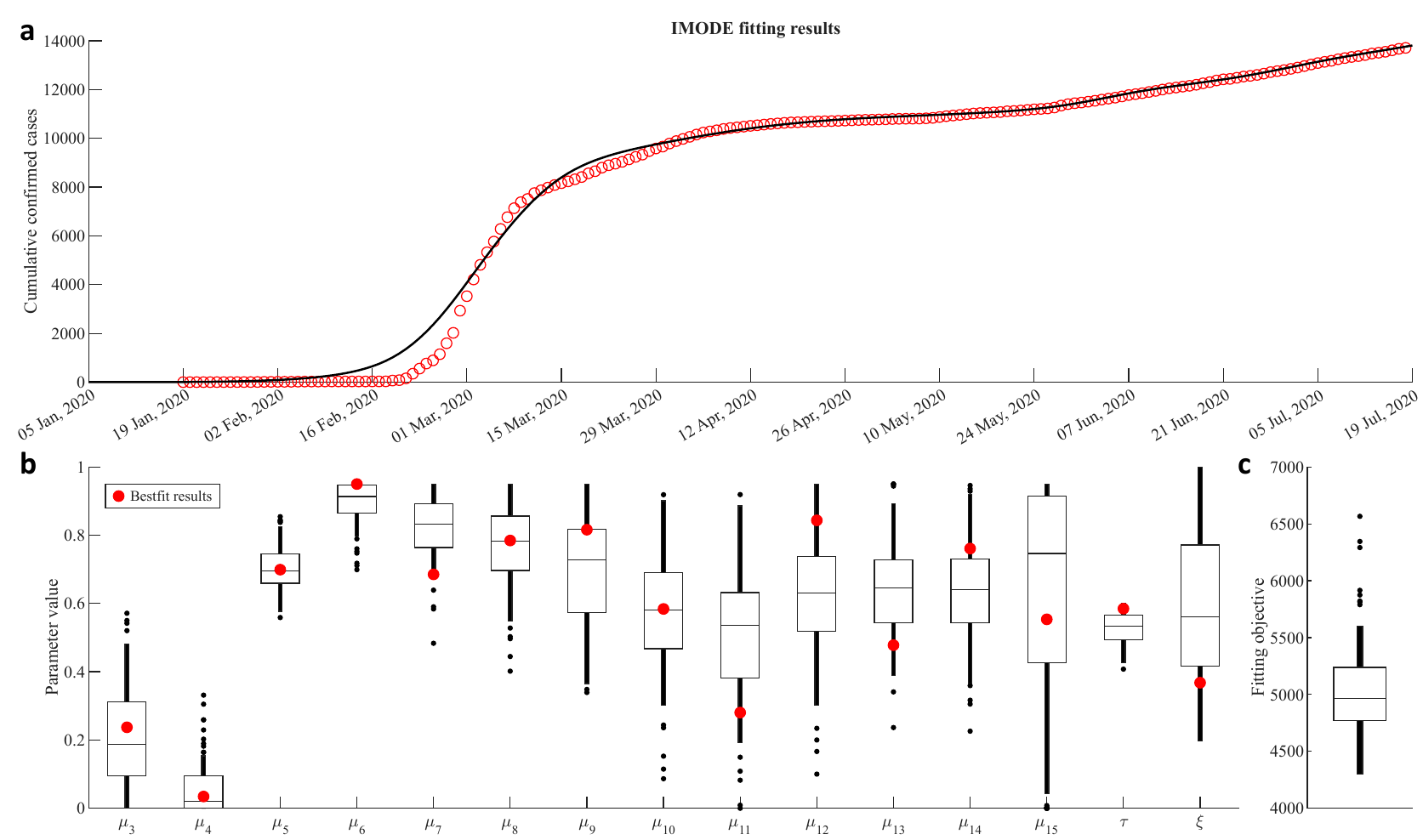}
    \caption{\textbf{IMODE estimation results.}
    (a) Cumulative confirmed cases and simulation results. (b) The 100 estimation results represented as a boxchart, with the best estimation results shown using red circle points. (c) L2-norm of the 100 estimation results.
    }\label{Afig01}
    \end{figure*}

    The IMODE captures the transmission dynamics of the epidemic, bridging the gap between theoretical modeling and observed data. As shown in the fitting results for cumulative confirmed cases, the model's best-fit trajectory closely aligns with the actual data recorded from January 2020 through July 2020, effectively reflecting the transition from the initial outbreak to the stable phases. This accuracy is supported by the optimization of a diverse set of policy-related parameters, where specific values for $\mu_{3},\,\mu_{4},\, \cdots,\,\mu_{14},\,\mu_{15}$, $\tau$, and $\xi$ were derived to minimize the discrepancy between the model and real-world observations. 
    
When a local optimizer is employed, user-specified initial values must often be considered, as the solution may be sensitive to the starting point. However, using a global optimizer eliminates this problem as it does not require an initial value for optimization. To reduce the randomness of the global optimizer, IMODE was run 100 times, and the best solution was selected \cite{IMODE}. Each simulation was terminated after 100,000 functional evaluations were performed. Figure~\ref{Afig01} shows the estimation results. Because the IMODE algorithm effectively searches for the most suitable solution within the loss landscape, its results were used as prior information in the MCMC process.

\subsection{MCMC process}\label{asec2_2}

The Markov chain Monte-Carlo (MCMC) algorithm was implemented using MATLAB \cite{MCMC,MCMC_package}. We set the prior distribution as a normal distribution, with the mean value obtained from the best IMODE results, to ensure that the Markov chains initiate within a high-probability region of the parameter space. A standard deviation of 0.05 was assigned to each policy-related parameter. This value was selected to provide an informative prior that stabilizes the search near the global optimum while maintaining enough variance to allow the delayed rejection adaptive metropolis (DRAM) algorithm \cite{MCMC_package} to accurately characterize the posterior distribution and its correlations. Unlike traditional non-informative priors, this hybrid approach utilizes IMODE-derived hyperparameters to refine the search space, which is particularly beneficial for high-dimensional policy parameters like $\mu(t)$ that vary over multiple two-week intervals.

\begin{figure*}[h]
    \centering
    \includegraphics[width=.95\textwidth]{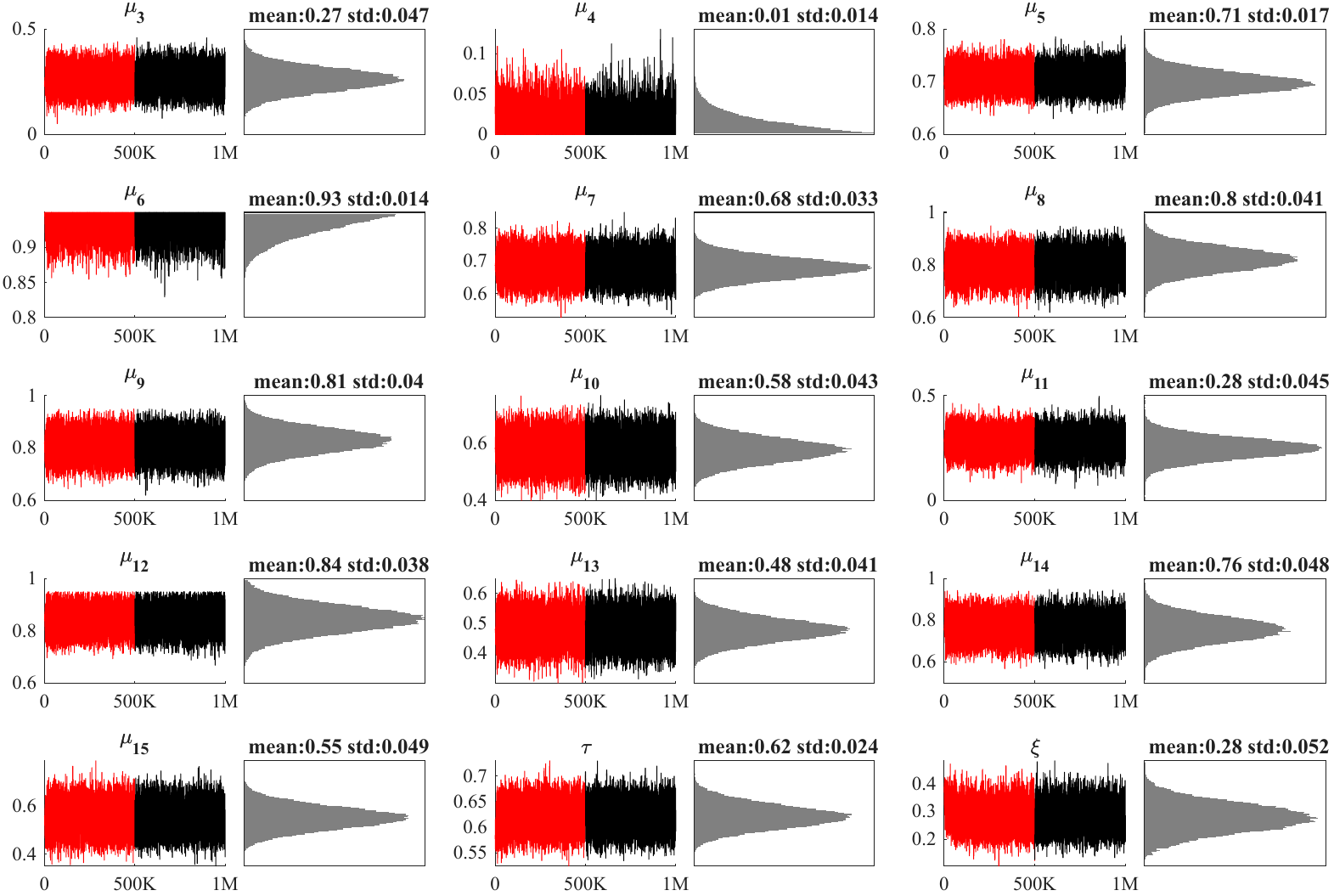}
    \caption{\textbf{MCMC chain and posterior distributions for epidemiological parameters.} 
    For each parameter, the chain illustrate the convergence of the chains over 1 million (1M) iterations. The red period indicates the burn-in, which constitutes half of the total iterations. The reported mean and standard deviation (std) values represent the posterior distribution of the estimated parameters.}\label{Afig03}
\end{figure*}

Figure~\ref{Afig03} illustrates the MCMC estimation results obtained using the DRAM algorithm. To quantify the uncertainty of the estimated parameters, we executed chains for 1 million iterations. The narrow standard deviations indicate high precision in the posterior distributions. Specifically, parameters associated with transmission and intervention efficacy showed distinct mean values, such as $\mu_{6}$ ($0.93 \pm 0.014$) and $\tau$ ($0.62 \pm 0.024$), providing a statistically grounded basis for simulating various pandemic scenarios. The high-probability regions identified by these chains, facilitated by the IMODE-derived priors, allowed for fast convergence to the posterior distribution.

\begin{figure*}[h]
    \centering
    \includegraphics[width=0.8\textwidth]{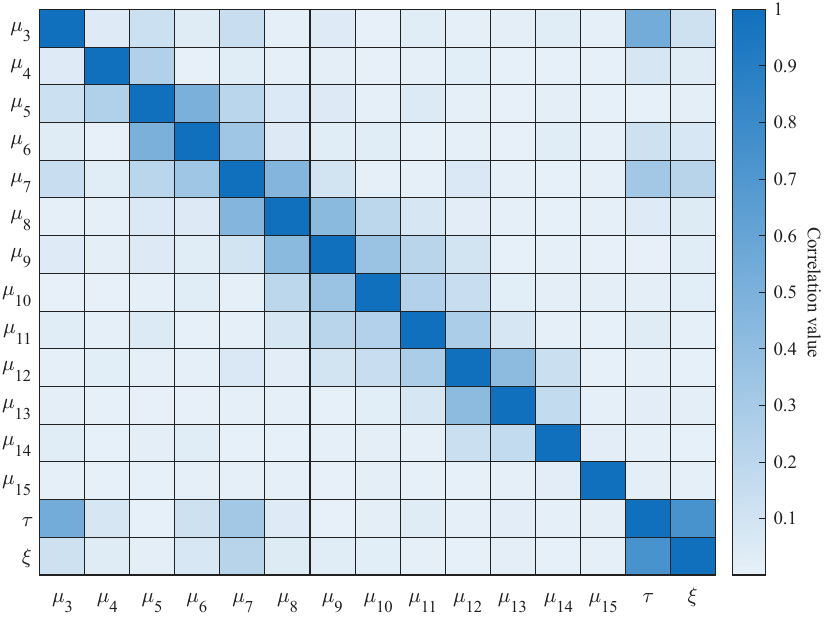}
    \caption{\textbf{Correlation between estimated parameters} 
    The heatmap illustrates the pairwise correlation coefficients between the posterior distributions of the model parameters, including $\mu_{3},\,\cdots,\,\mu_{15}$, $\tau$, and $\xi$.  
    }\label{Afig04}
\end{figure*}

To ensure the reliability of our estimation framework, we evaluated the correlations between the posterior distributions of all parameters generated by the MCMC process (Figure~\ref{Afig04}). The analysis reveals that the majority of parameters exhibit low to moderate correlations, with most correlation values remaining between 0 and 0.4. This indicates that the parameters are practically identifiable and that the IMODE-based prior selection effectively guided the MCMC chains toward distinct regions of the parameter space. The overall matrix confirms that the policy-related scaling factors, $\tau$ and $\xi$, are not strongly coupled with the transmission rate parameters.

\section{Sensitivity analysis results}\label{asec3}

    Sensitivity analysis results revealed the relative impact of each parameter on the model outputs, helping to identify which parameters can effectively control the outputs or which parameters may be negligible. We performed a partial rank coefficient correlation (PRCC) analysis of the number of infected and confirmed cases to assess the sensitivity of the model parameters. Notably, the confirmed cases were used as the outputs for parameter estimation. Among the policy-related parameters, $\mu$ and $\tau$ showed the highest sensitivity for both outputs and influenced the disease-related parameters $\beta$ and $\alpha$. 
    
\begin{figure*}[h]
    \centering
    \includegraphics[width=.95\textwidth]{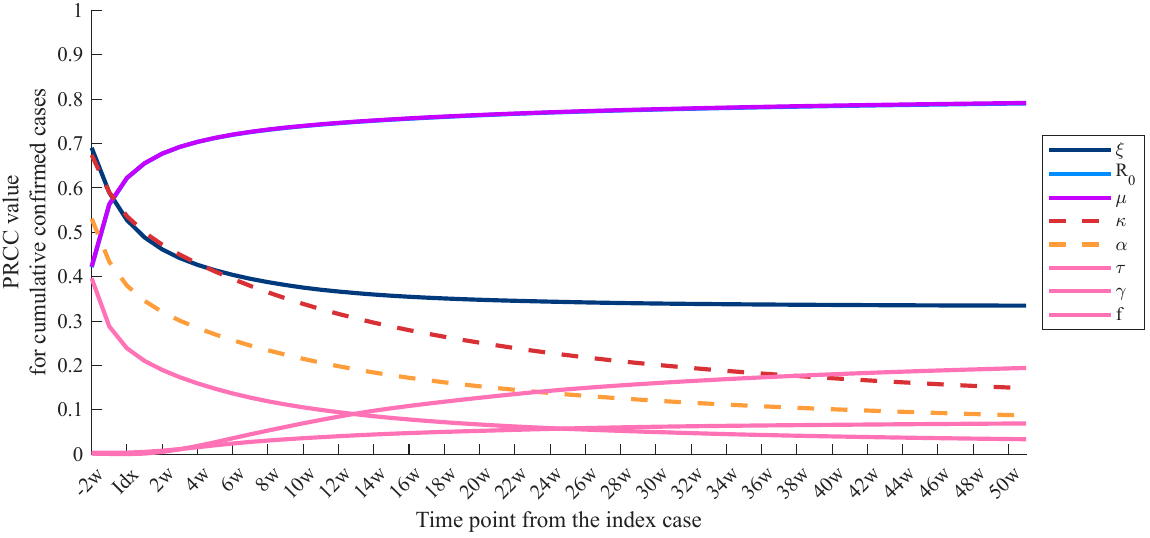}
    \caption{\textbf{Parameter estimation results.}
    Correlation between model simulation and parameters over time. Solid and dashed lines indicate positive and negative correlation, respectively.
    }\label{Afig05}
\end{figure*}
Sensitivity analysis assesses the impact of parameters on the simulation results. In this analysis, $\mu$ was set to a constant. Sensitivity was evaluated for two outputs: cumulative confirmed cases and number of infections. The cumulative cases were used for data-fitting, while the number of infections was employed for multi-objective optimization. Figure~\ref{Afig05} presents the PRCC results for these two outputs. $\beta$ has the highest correlation, exceeding 0.8 for all time. $\kappa$ and $\xi$ were 0.6723 and 0.6203, respectively, at the beginning of the simulation, but their values decreased monotonically to 0.2432 and 0.0837, respectively. The correlation for $\tau$ was 0.5 at the beginning of the simulation and increased to 0.8 by the end.
\end{appendices}

\newpage

\end{document}